\newcount\notenumber

\def\note{\advance\notenumber by 1
\footnote{$^{(\the\notenumber)}$}}

\def\E{{\cal E}}

\def\1{{\bf C}}

\def\X{{\bf X}}

\def\Q{{\bf Q}}

\def\N{{\bf N}}
\def\Z{{\bf Z}}

\def\x{{\bf x}}

\hyphenation{li-ne-a-ri con-si-de-ra di-stin-te Di-mo-stre-re-mo Sug-ge-ri-men-to
di-su-gua-glian-ze im-me-dia-ta-men-te i-nef-fet-ti-vo}

\hsize = 15truecm
\vsize = 22truecm

\hoffset = 0.6truecm
\voffset = 0.7truecm


\font\title=cmr10 scaled 1200


\centerline{\title On the length of the continued fraction for values of
quotients of power sums}\bigskip

\centerline{P. Corvaja\qquad U. Zannier}\bigskip

\noindent{\bf Abstract.} Generalizing a result of Pourchet, we show
that, if $\alpha,\beta$  are power sums over $\Q$ satisfying suitable
necessary assumptions,  the length of the continued fraction for
$\alpha(n)/\beta(n)$ tends to infinity as $n\rightarrow\infty$.
This will be derived from a uniform Thue-type inequality for
the rational approximations to the rational numbers $\alpha(n)/\beta(n)$,
$n\in\N$.\bigskip

\noindent{\bf Introduction.} The features of the continued fraction of
a positive real number  are usually extraordinarily difficult
to predict. However, if the numbers in question run 
through certain parametrized families, some regularity occasionally
appears. For instance, we may recall remarkable results by Schinzel
[S1], [S2], on the lengths of the periods of the continued fractions for
$\sqrt{f(n)}$, where
$f$ is a polynomial with rational coefficients and $n$ varies through
$\N$. 

In a similar direction, but much more simply, one can  deal with
the continued fractions of the {\it rational} numbers $r(n)$, $n\in\N$,
where now
$r\in\Q(X)$ is a rational function; for example, we may now look at the
length of the fraction: in [S2, Lemma 2], Schinzel proves that this
length is bounded as a function of $n$ (see also [M] for a more precise
result). Substantially, the reason for this is that 
$r(X)$ has a finite expansion as a {\it simple} continued fraction with
partial quotients in 
$\Q[X]$.\note{Observe however that such  continued fraction specialized
at
$X=n$ may not coincide with the continued fraction for $r(n)$.}

More generally, one may inquire about the length of the continued
fraction for the rational values of other   arithmetical
functions. In general, due to the lack of an Euclid's algorithm,
we have not a suitable ``functional continued  fraction"; this fact makes
the situation  different from the case of $\Q(X)$ and generally speaking 
considerably more difficult. 

For certain exponential functions,   we have a result by
Y. Pourchet, answering an appealing question of Mend\`es-France (see [M2, p.
214]); it states that, {\it for coprime  integers
$a,b>1$, the length of the continued fraction for
$(a/b)^n$ tends to infinity as $n\rightarrow\infty$}, in marked contrast
to the rational function case.

This theorem may be derived from the (deep)  lower bound  
$|(a/b)^n-(p/q)|\gg_\epsilon q^{-2}\exp(-\epsilon n)$  (any positive
$\epsilon$), where
$p,q$ are integers with $0<q<b^n$ (which follows e.g. from Ridout
generalization of Roth's Theorem [R]). The
estimate is in fact amply sufficient, since it implies that all the
partial quotients of the continued fraction in question are
$\ll_\epsilon
\exp(\epsilon n)$, and the statement about the length follows at once.
(See also [Z, Ex. II.6].)
\medskip 

In this direction, one may  consider  more general power sums in
place of $a^n,b^n$. Namely, in analogy to [CZ1] (where the notation
is slightly different) we consider the ring (actually a domain) 
$\E$ made up of the functions on $\N$ of the form
$$
\alpha(n)=c_1a_1^n+\ldots +c_ra_r^n,\eqno(1)
$$
where the coefficients $c_i$ are rational numbers, the $a_i$ are positive
integers  and the number $r$ of summands  is unrestricted.  (We shall
often normalize (1) so that the $c_i$ are nonzero and  the $a_i$ are
distinct; in this case the $a_i$ are called the ``roots" of $\alpha$.)

We then consider the ratio $\alpha/\beta$ of nonzero elements
of $\E$ and ask about the length of the continued fraction for the
values 
$\alpha(n)/\beta(n)$ (provided $\beta(n)\neq 0$, which is the case for
all  large
$n$).  It turns out
(Corollary 2 below) that the length of the relevant continued
fraction does not tend  to infinity if and only if $\alpha/\beta$ admits
an expansion as a (finite) continued fraction over
$\E$. Note that $\E$ is not euclidean, so in practice  this condition
``rarely"  holds. 

This result plainly generalizes Pourchet's, but   now the Ridout
Theorem seems no longer sufficient for a proof. In fact, we shall need
the full force of the  Schmidt Subspace Theorem, similarly to [CZ1],
where (among others) the related question of the integrality of the
values
$\alpha(n)/\beta(n)$ was investigated.

Similarly to the above sketched argument for Pourchet's Theorem, the
conclusions  about the continued fraction will be derived from a
``Thue-type" inequality for the rational approximations to the
values in question, which is the  content of the Theorem below.  Its
proof    follows   [CZ1] in
the main lines, except for  a  few technical points; however the
present  new simple  applications of those principles are perhaps not
entirely free of interest.

Preliminary to the statements, we introduce a little more 
notation.
\medskip

We let $\E_\Q$ be the domain of the functions of type
(1), but allowing the $a_i$ to be positive rationals. 
If $\alpha\in \E_\Q$ is written in the form (1), with nonzero $c_i$
and distinct $a_i$, we set $\ell(\alpha)=\max (a_1,\ldots ,a_r)$, agreeing
that $\ell(0)=0$.  

It is immediate to check that 
$\ell(\alpha\beta)=\ell(\alpha)\ell(\beta)$,
$\ell(\alpha+\beta)\le
\max(\ell(\alpha),\ell(\beta))$ and that $\ell(\alpha)^n\gg|\alpha(n)|\gg
\ell(\alpha)^n$ for $n\in\N$ tending to infinity (so in particular,
$\alpha(n)$ cannot vanish infinitely often if $\alpha\neq 0$).

\medskip
\noindent{\bf Theorem.}\ {\it Let $\alpha,\beta\in\E$ be nonzero and
assume that for all $\zeta\in\E$, $\ell(\alpha-\zeta\beta)\ge
\ell(\beta)$.   Then there exist    $k=k(\alpha,\beta)>0$ and
$Q=Q(\alpha,\beta)>1$ with the following properties. Fix
$\epsilon >0$; then,  for all but finitely many $n\in\N$ and 
for  integers
$p,q$,
$0<q<Q^n$,  we have 
$$
\left|{\alpha(n)\over\beta(n)}-{p\over q}\right|\ge 
{1\over q^k}\exp(-\epsilon n).
$$
}\medskip

\noindent We tacitly mean that the alluded finite set of exceptional
integers  includes those $n\in\N$ with $\beta(n)=0$;  this
finite set  may depend on
$\epsilon$.
\medskip

\noindent{\bf Remarks.} (i) The assumption about $\alpha,\beta$ expresses the lack of an
``Euclid division" in $\E$ for $\alpha:\beta$ and  cannot
be omitted. In
fact, suppose that
$\ell(\alpha-\zeta\beta)<\ell(\beta)$ for some $\zeta\in\E$. Then, the
values $\zeta(n)$ have bounded denominator  and verify
$|(\alpha(n)/\beta(n))-\zeta(n)|\ll \exp(-\epsilon_0n)$ for
large $n$ and some positive $\epsilon_0$ independent of $n$, against the
conclusion of the Theorem (with $p/q=p_n/q_n=\zeta(n)$).

Note that the assumption is automatic if $\ell(\alpha)\ge
\ell(\beta)$ and $\ell(\beta)$ does not divide
$\ell(\alpha)$.\smallskip  

(ii) By the short argument after  Lemma 1 below, given
$\alpha,\beta$,  one can test effectively whether a
$\zeta\in\E$ such that
$\ell(\alpha-\zeta\beta)<\ell(\beta)$ actually exist. Also, the proof
will show that a suitable $Q$ and ``exponent" $k$ may be computed. On the
contrary, the finitely many exceptional $n$'s cannot be computed with
the present method of proof.
\smallskip

(iii) Naturally, the lower bound  is not
significant   for  
$q$ larger than 
$\ell(\beta)^{n/(k-1)}$. So, some  upper
bound $Q^n$ for $q$ is not really restrictive. (Also, since the
dependence of
$k$ on
$\alpha,\beta$ is unspecified, and since 
$\ell(\beta)^{1/(k-1)}\rightarrow 1$ as
$k\rightarrow\infty$, it is immaterial here to specify
a suitable 
$Q$ in terms of $k$.)

Note also that   a lower bound $\gg Q^n$ for the
denominator of $\alpha(n)/\beta(n)$ follows, sharpening [CZ1, Thm. 1]; in
this direction, see [BCZ], [CZ3, Remark 1] and [Z, Thm. IV.3] for stronger
conclusions in certain special cases.
\smallskip

(iv) Like for the proofs in [CZ1], the method yields analogous
results for   functions of the form (1), but where $a_i$ are
any algebraic numbers  subject to the sole (but crucial)  restriction
that there exists a unique maximum among the absolute values $|a_i|$.  
Using the (somewhat complicated) method of [CZ2], one may relax this
condition, assuming only that not all the $|a_i|$ are equal. For the
sake of simplicity, we do not give the proofs of these results, which do
not involve new ideas compared to [CZ1], [CZ2] and the present paper, but
only complication of detail.\smallskip

(v) The Theorem can be seen as a  Thue-type inequality with
{\it moving targets} (similarly to [C] or [V]). It seems an interesting,
but difficult, problem, to obtain (under suitable necessary
assumptions) the ``Roth's exponent"
$k=2$, or even some exponent  independent of $\alpha,\beta$. (As
recalled above, this holds in the special cases of Pourchet's Theorem.)
\medskip

\noindent{\bf Corollary 1.} \ {\it Let $\alpha,\beta$ be as in the
Theorem.  Then the length of the continued fraction for
$\alpha(n)/\beta(n)$ tends to infinity as
$n\rightarrow\infty$.}\medskip

This corollary is in fact a lemma for the following result, which gives a
more precise description.\medskip

\noindent {\bf Corollary 2}. {\it Let $\alpha,\beta\in \E$ be nonzero.
Then the length of the continued fraction for
$\alpha(n)/\beta(n)$ is bounded for infinitely many $n\in\N$ if and only
if   there exist  power sums $\zeta_0,\ldots,\zeta_{k}\in\E$
such that we have the identical continued fraction expansion
$$
{\alpha(n)\over\beta(n)}=[\zeta_0(n),\zeta_1(n),\ldots,\zeta_{k}(n)].
$$
If this is the case, the mentioned length is uniformly bounded for all
$n\in\N$. }
\medskip

It will be pointed out how the  condition on $\alpha/\beta$ may be
checked effectively.

The special case $\alpha(n)=a^n-1$, $\beta(n)=b^n-1$ appears as [Z, Ex.
IV.12]. For completeness we give here this application.\medskip

\noindent {\bf Corollary 3}. {\it Let $a,b$ be multiplicatively
independent positive integers. Then the length of the continued fraction
for  $(a^n-1)/(b^n-1)$
tends to infinity with $n$.
}\bigskip

\noindent{\bf Proofs.} We start with the
following very simple:
\medskip

\noindent{\bf Lemma 1.}\ {\it Let $\alpha,\beta\in\E$ be nonzero and let
$t$ be any positive number. Then there exists  $\eta\in\E_\Q$ 
 such that $\ell(\alpha-\eta\beta)<t$. Such an $\eta$ may be computed in
terms of $\alpha,\beta,t$.}\medskip

\noindent{\it Proof.} Write $\beta(n)=cb^n(1-\delta(n))$, where
$c\in\Q^*$, $b=\ell(\beta)$ and where $\delta\in\E_\Q$ satisfies
$u:=\ell(\delta)<1$. In particular, we have $|\delta(n)|\ll u^n$, so for
a fixed integer $R$ we have an approximation, for $n\rightarrow\infty$, 
$$
{\alpha(n)\over\beta(n)}={\alpha(n)\over cb^n}
\left(\sum_{r=0}^R\delta(n)^r\right)+
O\left(\left({\ell(\alpha)u^{R+1}\over b}\right)^n\right).
$$
Choose $R$ so that $u^R<t/u\ell(\alpha)$ and define $\eta$ by the first
term on the right. Then $\eta\in\E_\Q$ and 
$|\alpha(n)-\eta(n)\beta(n)|\ll (\ell(\alpha)u^{R+1})^n$. Therefore
$\ell(\alpha-\eta\beta)\le \ell(\alpha)u^{R+1}<t$,  concluding the
proof.\medskip

We note at once that this argument may be used to check the
assumption for the Theorem. In fact, suppose
$\ell(\alpha-\zeta\beta)<\ell(\beta)$ for a $\zeta\in\E$. Then, if $\eta$
is as in the lemma, with $t=\ell(\beta)$,   it follows that
$\ell(\eta-\zeta)<1$. Since $\eta$ may be constructed and since the
``roots" of $\zeta$ are positive integers by assumption, this
inequality determines $\zeta$ uniquely; namely, if a suitable
$\zeta$ exists, it is just the ``subsum" of $\eta$ made up  with the
integer roots.\medskip

To go on, for the reader's convenience, we recall a version of Schmidt's
Subspace Theorem suitable for us; it is  due to H.P. Schlickewei (see
[Schm2, Thm.1E, p.178]).\medskip

\noindent{\bf Subspace Theorem.}Ê\ {\it Let $S$ be a finite set of
places of $\bf Q$, including the infinite one and normalized in
the usual  way (i.e.  $|p|_v=p^{-1}$ if $v|p$). 
For $v\in S$ let $L_{0v},\ldots ,L_{hv}$ be $h+1$ linearly independent
linear forms in $h+1$ variables with rational coefficients and let
$\delta >0$.  Then the
solutions ${\bf x}:=(x_0,\ldots ,x_h)\in{\bf Z}^{h+1}$ to the inequality
$$
\prod_{v\in S}\prod_{i=0}^h|L_{iv}({\bf x})|_v\le ||{\bf x}||^{-\delta}
$$
where $||{\bf x}||:= \max\{|x_i|\}$, are all contained in finitely many
proper  subspaces of ${\bf Q}^{h+1}$.}\medskip

For our last lemma we could invoke  results by Evertse 
[E]; for completeness we give a short  proof of the special case we
need.\medskip

\noindent{\bf Lemma 2.} \ {\it Let $\zeta\in\E_\Q$ and let $D$ be the
minimal positive integer such that $D^n\zeta(n)\in\E$. Then, for every
$\epsilon >0$ there are only finitely many $n\in\N$ such that the
denominator of $\zeta(n)$ is smaller than $D^n\exp(-\epsilon
n)$.}\medskip

\noindent{\it Proof.} We may write
$D^n\zeta(n)=\sum_{i=0}^hu_ie_i^n$, where $u_i\in\Q^*$ and $e_i$ are
distinct positive integers with $(e_0,\ldots ,e_h,D)=1$.  Let $S$ be the
set of places of $\Q$ made up of the infinite one and of those dividing
$De_0\cdots e_h$. 

Define the linear form $L(\X):=\sum_{i=0}^hu_iX_i$.

For each $v\in S$ we define linear forms $L_{iv}$, $i=0,\ldots ,h$, as
follows. If $v$ does not divide $D$ (including $v=\infty$), then we put
simply
$L_{iv}=X_i$ for
$i=0,\ldots ,h$. If $v|D$, then there exists $j=j_v$ such that
$|e_j|_v=1$; then we set $L_{jv}=L$ and $L_{iv}=X_i$ for $i\neq j$. The
forms $L_{iv}$, $i=0,\ldots ,h$, are plainly linearly independent for
each fixed
$v\in S$.

Suppose now that the conclusion does not hold, so for some positive
$\epsilon$ and for all
$n$ in an infinite set $\Sigma\subset\N$ the denominator of $\zeta(n)$ is
$\le D^n\exp(-\epsilon n)$. For $n\in\Sigma$, put 
$\x=\x(n)=(e_0^n,\ldots ,e_h^n)$.

Then, for $n\in\Sigma$, the numerator  of
$D^n\zeta(n)$ has a g.c.d. with $D^n$ which is $\ge\exp(\epsilon n)$. In
turn this gives 
$$
\prod_{v|D}|L(\x(n))|_v\le\exp(-\epsilon n),\qquad n\in\Sigma.
$$
Also, for $n\in\Sigma$, we have
$$
\prod_{v\in S}\prod_{i=0}^h|L_{iv}(\x)|_v=
\left(\prod_{v|D}|L(\x)|_v\right)\prod_{v\in S}\prod_{i=0}^h|e_i^n|_v
=\prod_{v|D}|L(\x)|_v,
$$
where the first equality holds because   $|e_{j_v}|_v=1$ for $v|D$
and the second because the (inner) double product is
$1$ by the product formula. 
Hence
$$
\prod_{v\in S}\prod_{i=0}^h|L_{iv}(\x)|_v\le\exp(-\epsilon n),\qquad
n\in\Sigma.
$$
 Since however
$||\x(n)||$ is bounded exponentially in $n$, we may apply the Subspace
Theorem (with a suitable positive $\delta$) and conclude that all the
$\x(n)$,
$n\in\Sigma$, lie in a finite union of proper subspaces. But
$e_0,\ldots ,e_h$ are distinct positive integers, whence each subspace
corresponds to at most finitely many $n$, a contradiction which
concludes the proof.

\medskip

\noindent{\it Proof of Theorem.} Let $\eta$ be as in Lemma 1, with
$t=1/9$. Write
$$
\eta(n)=g_1({e_1\over d})^n+\ldots +g_h({e_h\over d})^n,
$$
where $g_i\in\Q^*$ and where $e_i,d$ are positive integers with
$e_1>e_2>\ldots >e_h$. 

We define $k:=h+3$, $Q=\exp(1/k)$, and we assume
that  for some fixed
$\epsilon >0$ (which we  may take $<1/6$, say) there exist infinitely
many triples
$(n,p,q)$ of integers with $0<q\le Q^n$, $n\rightarrow\infty$ and 
$$
\left|{\alpha(n)\over\beta(n)}-{p\over q}\right|\le  
{1\over q^k}\exp(-\epsilon n).\eqno(2)
$$
We
shall eventually obtain a contradiction, which will prove what we
want.\medskip

We proceed to define the data for an application of the Subspace
Theorem. We let $S$ be the finite set of    places of $\Q$ consisting of
the infinite one and of all the places dividing $de_1\cdots e_h$. 

We  define   linear forms in $X_0,\ldots ,X_h$ as follows. For
$v\neq\infty$ or for $i\neq 0$ we set simply $L_{iv}=X_i$. We define  the
remaining form $L_{0\infty}$ by 
$$
L_{0\infty}(\X)=X_0-g_1X_1-\ldots -g_hX_h.
$$
Plainly, for each $v$ the forms $L_{iv}$ are independent. For a triple
$(n,p,q)$ as above we set 
$$
\x=\x(n,p,q)=(pd^n,qe_1^n,\ldots ,qe_h^n)\in\Z^{h+1}
$$
and we proceed to estimate the double product $\prod_{v\in
S}\prod_{i=0}^h|L_{iv}(\x)|_v$. 

For $i\neq 0$ we have 
$\prod_{v\in S}|L_{iv}(\x)|_v=\prod_{v\in S}|qe_i^n|_v\le q$. In fact,
we have $\prod_{v\in S}|e_i^n|_v=1$ by the product formula, since
the $e_i$ are $S$-units; also, $\prod_{v\in S}|q|_v\le |q|$, since $q$ is
an integer.

Further, we have $\prod_{v\in
S}|L_{0v}(\x)|_v={|L_{0\infty}(\x)|\over|x_0|}\prod_{v\in S}|x_0|_v$.
As before, we see that the inner product is bounded by $|p|$, since $d^n$
is an
$S$-unit.

Also,   $|L_{0\infty}(\x)|= qd^n|\eta(n)-(p/q)|$. Then, using (2) and
  $|(\alpha(n)/\beta(n))-\eta(n)|\ll (t/\ell(\beta))^n\le t^n$, we obtain
$$
|L_{0\infty}(\x)|\ll qd^nt^n+qd^n(q^{-k}\exp(-\epsilon n)).
$$
Since $q^k\le Q^{kn}=\exp(n)$ and since  $t=1/9<\exp(-1-\epsilon)$, the
first term on the right does not exceed the second one, whence
$$
|L_{0\infty}(\x)|\ll qd^n(q^{-k}\exp(-\epsilon n)).
$$
Combining this with the previous bound we get
$$
\prod_{v\in S}|L_{0v}(\x)|_v\ll qd^n(q^{-k}\exp(-\epsilon
n))|x_0|^{-1}|p|
= q^{-(k-1)}\exp(-\epsilon n).
$$
Finally, for the double product this yields (since $k>h+1$)
$$
\prod_{v\in S}\prod_{i=0}^h|L_{iv}(\x)|_v
\ll q^{-(k-h-1)}\exp(-\epsilon n)\ll \exp(-\epsilon n).
$$
Now, (2) entails $|p|\ll q\ell(\alpha)^n\ll (Q\ell(\alpha))^n$, whence 
$||\x||\ll C^n$ for some $C$ independent of $n$. Hence the double
product is bounded by $||\x||^{-\delta}$ for some fixed positive
$\delta$ and large enough $n$. 

The Subspace Theorem then implies that all the vectors $\x$ in question
are contained in a certain finite union of proper subspaces of
$\Q^{h+1}$. 

In particular, there exists a fixed subspace, say of
equation $z_0X_0-z_1X_1-\ldots -z_hX_h=0$, containing an infinity of the
vectors in question.  We cannot have $z_0=0$, since this would entail 
$\sum_{i=1}^hz_ie_i^n=0$ for an infinity of $n$; in turn, the fact that
the
$e_i$ are distinct would imply $z_i=0$ for all $i$, a contradiction.
Therefore we may assume that $z_0=1$, and we find that, for the
$n$'s corresponding to the vectors in question,
$$
{p\over q}=\sum_{i=1}^hz_i({e_i\over d})^n=\zeta(n),\eqno(3)
$$
say, where $\zeta\in\E_\Q$. We now show that actually $\zeta$   lies in
$\E$. Assume the contrary; then the minimal positive integer $D$ so that
$D^n\zeta(n)\in\E$ is $\ge 2$. But then equation (3) together with Lemma
2 implies that $q\gg 2^n\exp(-\epsilon n)$. Since this would hold for
infinitely many $n$, we would find
$Q\ge  2\exp(-\epsilon)$, whence 
$\exp((1/k)+\epsilon)\ge 2$, a contradiction. (Recall that $k\ge
3,\epsilon <1/6$.)

Therefore $\zeta$ lies in $\E$; using (3) to substitute in (2) for $p/q$
we find that, for  an infinity of $n$,
$$
\left|{\alpha(n)\over \beta(n)}-\zeta(n)\right|
\le \exp(-\epsilon n).
$$
In particular, $\ell((\alpha/\beta)-\zeta)<1$, whence
$\ell(\alpha-\zeta\beta)<\ell(\beta)$, contrary to the
assumptions,   concluding the proof.\bigskip

\noindent{\it Proof of Corollary 1.} For notation and basic facts about
continued fractions we refer to   [Schm1, Ch. I] and let
$p_r(n)/q_r(n)$,
$r=0,1,\ldots$,  be the (finite) sequence of convergents of the continued
fraction for
$\alpha(n)/\beta(n)$, where we may assume
that $\alpha(n)$ and
$\beta(n)$ are positive. 

 As is well-known,
$p_r(n),q_r(n)$ are positive integers, $q_0(n)=1$,  and we
have 
$$
\left|{\alpha(n)\over\beta(n)}-{p_r(n)\over q_r(n)}\right|\le 
{1\over q_{r+1}(n)q_r(n)}.\eqno(4)
$$
Suppose the conclusion false, so for some fixed $R\in \N$ there is an
infinite set
$\Sigma\subset\N$ such that for $n\in\Sigma$ the continued fraction for 
$\alpha(n)/\beta(n)$ has finite length $R+1\ge 1$. By this we mean that
$p_R(n)/q_R(n)=\alpha(n)/\beta(n)$.

Since
$\alpha,\beta$ satisfy the assumptions for the Theorem, there exist  
$Q>1$, $k\ge 2$  as in that statement. 

Define now the sequence
$c_0,c_1,\ldots
$ by $c_0=0$ and $c_{r+1}=(k-1)c_r+1$ and choose  
 a positive number $\epsilon <c_R^{-1}\log Q$, so 
$\exp(c_R\epsilon) <Q$;   we proceed to show by induction on
$r=0,\ldots , R$, that, for large $n\in\Sigma$,
$$
q_{r}(n)\le \exp(c_r\epsilon n).\eqno(5)
$$
We have $q_0(n)=1$, so (5) is true  for $r=0$; we shall now show that 
(if $r\le R-1$), for large $n$,  (5) implies  the same inequality with
$r+1$ in place of
$r$.

Observe first that, by construction,  $\exp(c_r\epsilon)\le
\exp(c_R\epsilon)<Q$, whence $q_r(n)\le \exp(c_r\epsilon n)<Q^n$ for
large enough
$n\in\Sigma$. We may
then apply the Theorem with $p=p_r(n)$, $q=q_r(n)$ and deduce that  for
large
$n\in\Sigma$ we have 
$$
\left|{\alpha(n)\over\beta(n)}-{p_r(n)\over q_r(n)}\right|\ge
q_r(n)^{-k}\exp(-\epsilon n).
$$
Hence, combining with (4),
$$
q_{r+1}(n)\le q_r(n)^{k-1}\exp(\epsilon n)\le\exp(((k-1)c_r+1)\epsilon n)
=\exp(c_{r+1}\epsilon n)
$$
by the inductive assumption and the definition
of $c_{r+1}$. This induction proves (5) for $r\le R$.

Finally, by (5) with $r=R$, we have $q_R(n)\le\exp(c_R\epsilon n)$ for
large $n\in\Sigma$, whence  $q_R(n)<Q^n$,
since 
$\exp(c_R\epsilon )<Q$ by construction. But then the Theorem holds for
$p=p_R(n)$, $q=q_R(n)$, leading in particular to
$0=|(\alpha(n)/\beta(n))-(p_R(n)/q_R(n))|>0$, a contradiction which
concludes the argument.\medskip

\noindent{\it Proof of Corollary 2.} For convenience let us denote by
$\psi (x)$ the length of the continued fraction for the rational
number $x$.

We start with the hardest half of the proof; namely, assuming that
$\psi(\alpha(n)/\beta(n))$ is bounded for $n$ in an infinite sequence
$\Sigma\subset\N$, we  prove that $\alpha/\beta$ has a finite
continued fraction over $\E$.

Let us argue by induction on $\ell(\alpha)+\ell(\beta)$, the result
being trivial when this number is $\le 2$.  

Since $|\psi(x)-\psi(1/x)|\le 1$, we may plainly assume $\ell(\alpha)\ge
\ell(\beta)$. By our present assumption, the hypothesis of   Corollary 
1 cannot hold for
$\alpha, \beta$, so there exists $\zeta\in\E$ such that, putting
$\eta=\alpha-\zeta\beta$, we have 
$\ell(\eta)<\ell(\beta)$. For large $n$, we have $\beta(n)\neq 0$ and 
$$
{\eta(n)\over\beta(n)}={\alpha(n)\over \beta(n)}-\zeta(n).
$$
Since $\zeta\in\E$, the rational numbers $\zeta(n)$, $n\in\N$,  have a
finite common denominator $D> 0$. We now appeal to [S1, Thm. 1], which
bounds $\psi((p/q)-r)$ in terms of $\psi(p/q)$, where $r$ is a rational
number with fixed denominator. 
Applying this with $p/q=\alpha(n)/\beta(n)$, $r=\zeta(n)$ and noting
that  $\psi(\alpha(n)/\beta(n))$ is bounded by assumption for
$n\in\Sigma$,  we immediately obtain that the numbers
$\psi(\eta(n)/\beta(n))$ are also bounded for
$n\in\Sigma$. Since
$\ell(\eta)+\ell(\beta)<\ell(\alpha)+\ell(\beta)$, the inductive
assumption implies the existence of a finite continued fraction
expansion over $\E$ for $\eta/\beta$. In turn, this yields the same for
$\alpha/\beta$, proving the conclusion.

The ``converse" part of Corollary 2 is again immediate from [S1, Thm.
1].\medskip

Like  for Corollary 1, the condition that $\alpha/\beta$ admits a finite
continued fraction over $\E$ may be checked effectively. 
To verify this,
observe first that, if $\alpha/\beta$ admits a continued fraction
over $\E$ as in the statement, then it also admits an expansion such that
$\ell(\zeta_i)\ge 1$ for all $i\ge 1$; this claim follows at once from
the  
  continued fraction identity
$$
[A,c,B]=[A+c^{-1},-c(cB+1)],
$$
which allows to absorb all possible ``constant" partial quotients $c$.

With this proviso (but not otherwise) the continued fraction over $\E$ is
unique (if it exists); its existence may be checked e.g. by iterating the
effective criterion for checking the assumption of the Theorem (or of 
Corollary 1), as in the comments after Lemma 1.\medskip

\noindent{\it Proof of Corollary 3.} We could invoke Corollary 2,
but it is perhaps simpler to argue directly, by induction on $a+b$ (the
assertion being empty for $a+b=2$). On exchanging $a,b$ we may assume
that $a>b>1$. 

Now, we may assume that $a$ is divisible by $b$; in fact,
if this does not hold, we plainly have 
$\ell(a^n-1-\zeta(b^n-1))\ge a>b$ for any
$\zeta\in\E$; in turn,  we may then apply  ({\it a fortiori}) Corollary 1
to  $\alpha(n)=a^n-1$, $\beta(n)=b^n-1$, yielding the present
assertion.

Write then $a=bc$, where   $c$ is a positive integer $<a$.
Note that $b,c$ are multiplicatively independent (since $a,b$ are such)
whence the continued fraction for $(c^n-1)/(b^n-1)$ has length tending
to infinity, by the inductive assumption. But then the identity
$$
{a^n-1\over b^n-1}=c^n+{c^n-1\over b^n-1} 
$$
immediately implies the sought conclusion.\medskip

\bigskip

\noindent{\bf Final remarks.} (a) It seems a very difficult problem to
quantify the corollaries, namely to prove an explicit lower bound for the
length of the relevant continued fractions. Even in the special case of
Pourchet's Theorem, this seems to be related with an explicit version of
Ridout's (or Roth's) Theorem, which presently appears inaccessible.
\smallskip

(b) The methods of proof probably lead  to results similar to the Theorem,
where $\alpha(n)/\beta(n)$ is replaced by some  algebraic function of a
power sum, like e.g. the $d$-th root $\alpha(n)^{1/d}$. The pattern
should combine the arguments from, say, [CZ1, Cor. 1] (or [CZ4, Thm. 2])
with the above ones. 

In case $d=2$ it is perhaps possible to obtain   corollaries similar to
the present one, but in the direction of Schinzel's mentioned papers;
namely, under suitable assumptions on the power-sum $\alpha\in\E_\Q$,
{\it the length of the period of the continued fraction for
$\sqrt{\alpha(n)}$  tends to infinity with
$n$}.

\bigskip

\noindent{\bf References.}\medskip

\item{[BCZ]} -  Y. Bugeaud, P. Corvaja, U. Zannier, An upper bound for
the G.C.D. of $a^n-1$ and $b^n-1$, {\it Math. Z.} {\bf 243} (2003),
79-84.\smallskip 

\item{[C]} - P. Corvaja, Une application nouvelle de la
m\'ethode de Thue, {\it Annales de l'Institut Fourier}
{\bf 45}, (1995), 1177-1203.\smallskip

\item{[CZ1]} - P. Corvaja, U. Zannier,  Diophantine  equations with
power sums and Universal Hilbert Sets, {\it Indag. Mathem., N.S.,} {\bf
9} (3) (1998), 317-332.\smallskip

\item{[CZ2]} - P. Corvaja, U. Zannier, Finiteness of integral values for
the ratio of two linear recurrences, {\it Invent. Math.} {\bf 149}
(2002), 431-451.\smallskip

\item{[CZ3]} - P. Corvaja, U. Zannier, On the greatest prime factor of
$(ab+1)(ac+1)$,   {\it Proc. Amer. Math. Soc.}, {\bf 131} (2003),
1705-1709.\smallskip

\item{[CZ4]} - P. Corvaja, U. Zannier, Some New Applications of the
Subspace Theorem, {\it Compositio Math.} {\bf 131} (2002),
319-340.\smallskip

\item{[E]} - Y.-H. Evertse, On sums of $S$-units and linear
recurrences, {\it Compositio Math.}, {\bf 53} (1984), 225-244.\smallskip

\item{[M]} - M. Mend\`es-France, Quelques probl\`emes relatifs \`a
la th\'eorie des fractions continues limit\'ees, {\it S\'em.
Delange-Pisot-Poitou, 13e ann\'ee (1971/72), Th\'eorie des nombres, 1,
Exp. No. 2, Secr\'etariat Math\'ematique, Paris}, 1973.\smallskip

\item{[M2]} - M. Mend\`es-France, Sur les fractions continues
limit\'ees, {\it Acta Arith.} {\bf 23} (1973), 207-215.\smallskip

\item{[R]} - D. Ridout, The $p$-adic generalization of the
Thue-Siegel-Roth theorem, {\it Mathematika}, {\bf 5} (1958),
40-48.\smallskip

\item{[S1]} - A. Schinzel, On some problems of the
arithmetical theory of continued fractions, {\it Acta Arith.}
{\bf 6} (1961), 393-413.\smallskip

\item{[S2]} - A. Schinzel, On some problems of the
arithmetical theory of continued fractions II, {\it Acta
Arith.} {\bf 7} (1962), 287-298.\smallskip

\item{[Schm1]} - W.M. Schmidt, {\it Diophantine Approximation}, Springer
LNM 785, 1980.\smallskip

\item{[Schm2]} - W.M. Schmidt, {\it Diophantine Approximations and 
Diophantine
 Equations}, Springer Verlag  LN {\bf 1467}, 1991.\smallskip

\item{[V]} - P. Vojta, Roth's theorem with moving targets, {\it
Internat. Math. Res. Notices}, {\bf 3} (1996)
109-114.\smallskip

\item{[Z]} - U. Zannier, {\it Some applications of Diophantine
Approximation to Diophantine Equations}, to appear in the series
``Dottorato di Ricerca in Matematica", ETS, Pisa.\smallskip

\bigskip

\vfill

Pietro Corvaja

Dipartimento di Matematica, Universit\`a di Udine

via delle Scienze 206, 33100 Udine (Italy)

corvaja@dimi.uniud.it\medskip

Umberto Zannier

Istituto Univ. di Architettura, D.C.A.

S. Croce 191, 30135 Venezia (Italy)

zannier@brezza.iuav.unive.it

\end